\documentclass{article}
\usepackage[T2A]{fontenc}
\usepackage[cp1251]{inputenc}
\usepackage[english,russian]{babel}
\usepackage[tbtags]{amsmath}
\usepackage{amsfonts,amssymb,mathrsfs,amscd}
\usepackage[hyper]{msb-a}
\JournalName{}
\numberwithin{equation}{section}

\begin{document}

\title{Topology of isoenergy surfaces of Kovalevskaya integrable case on the Lie algebra so(4)}
\author
{Vladislav~Kibkalo}
\address{Faculty of Mechanics and Mathematics of Lomonosov Moscow State University}
\email{slava.kibkalo@gmail.com}

\udk{517.938.5}

\maketitle

\begin{fulltext}

\begin{abstract}
In the paper we determine the class of diffeomorphism of three-dimensional regular common level surfaces of Hamiltonian and Casimir functions for the analog of Kovalevskaya case on Lie algebra $\textrm{so}(4)$. We start from Fomenko-Zieschang invariants of Lioville foliations on these manifolds that were calculated by the author earlier.

The bibliography: 8 items.
\end{abstract}

\markright{Topology of isoenergy surfaces of Kovalevskaya case so(4)}

\section{Kovalevskaya case and its analogues}

Classical Kovalevskaya top was discovered and described by S. Kovalevskaya [1] and [2]. This dynamical system is integrable case of the Euler equations on the Lie algebra $\operatorname{e}(3)$. It was included by I.Komarov [3] in a one-parameter family of integrable systems on the pencil of Lie algebras $\operatorname{so}(4)-\operatorname{e}(3)-\operatorname{so}(3, 1)$.

The Lie--Poisson brackets is the following \begin{equation}  \{J_i, J_j\} = \varepsilon_{ijk}J_k, \quad \{J_i, x_j\} = \varepsilon_{ijk}x_k,
\quad \{x_i, x_j\} = \varkappa \varepsilon_{ijk}J_k, \end{equation}
where $\varepsilon_{ijk}$ is symbol of permutation and $\varkappa \in \mathbb{R}$. The cases $\varkappa>0, \varkappa =0$ and $\varkappa <0$ correspond to the Lie algebras $\operatorname{so}(4), \operatorname{e}(3)$ and $\operatorname{so}(3, 1)$ respectively.

Every regular 4-dimensional common level surface  $M^4_{a, b} = \left\{ f_1 =a, f_2 = b\right\} $ of the Casimir functions of the brackets (1.1) is a symplectic leaf and is diffeomorphic to $S^2 \times S^2$ in the case of $\varkappa >0$.

Bracket (1.1) has the following Casimir functions $f_1$ and $f_2$ for $\varkappa \in \mathbb{R}$:
\begin{equation}  f_1 = (x_1^2 + x_2^2 + x_3^2) + \varkappa (J_1^2 +J_2^2 +J_3^2),
\qquad  f_2 = x_1 J_1 + x_2 J_2 +x_3 J_3.
\end{equation}

Every their level surface $M^4_{a, b} = \left\{ f_1 =a, f_2 = b\right\}$ is four-dimensional and regular in the case of $a > \sqrt{\varkappa} |b|$ for the case of positive $\varkappa$, $a >0$ for the case of $\varkappa =0$ and for all $(a, b)$ excepting $(0, 0)$ for the case of $\varkappa <0$. Such $M^4$ are diffeomorphic to $S^2 \times S^2$ if $\varkappa >0$ and $T^{*} S^2$ otherwise.

The following smooth function $H$ define an integrable Hamiltonian system with 2 degrees of freedom on every regular $M^4_{a, b}$ and arbitrary $\varkappa \in \mathbb{R}$ w.r.t. the bracket (1.1) with first integral $K$
\begin{equation}
H = J_1^2 + J_2^2 + 2J_3^2 + 2 c_1 x_1,
\end{equation}
\begin{equation} K = (J_1^2 - J_2^2-2c_1 x_1 + \varkappa
c_1^2)^2 + (2J_1 J_2 - 2 c_1 x_2)^2.
\end{equation} Here $c_1$ is an arbitrary constant. We can assume that $c_1 = 1$ and $\varkappa \in \{-1, 0, 1\}$ (see e.g. [4]).

Common level surface of Casimir functions (1.2) and Hamiltonian (1.3) of a system is called an isoenergy surface. We say that such submanifold $Q^3_{a, b, h}$ is \textit{regular} if
\[Q^3_{a, b, h} = \{x \in M^4_{a, b} | H(x) = h\}\]
\begin{itemize}
\item it is three-dimensional and $\mathrm{d} H \ne 0$ on $Q^3_{a, b, h}$,
\item every $x \in Q^3_{a, b, h}$ is not of the rank $0$ for the momentum map $(H, K)$,
\item every rank 1 critical point is of the Bott type (i.e. nondegenerate).
\end{itemize}

In this paper we will calculate the class of diffeomorhism of isoenergy manifolds for Kovalevskaya case on Lie algebra $\operatorname{so}(4)$ can be calculated starting from their Fomenko-Zieschang invariants. All the necessary information about the Fomenko theory on the topological analysis of integrable Hamiltonian systems used in this paper can be found in [5].

In [6] Fomenko-Zieschang invariants for the classical Kovalevskays case were calculated by A.Bolsinov, A. Fomenko and P.Richter using the loop molecule approach for expression of some bases through uniquely defined cycles. In [7] by similar way V. Kibkalo expressed admissible coordinate systems that appear only in the case of $\operatorname{so}(4)$.

In [8] V.kibkalo calculated Fomenko-Zieshcang invariants for Liouville foliation on every regular isoenery manifold $Q^3_{a, b, h}$ and described 3-dimensional connected sets of triples $(a, b, h)$ in the space $\mathbb{R}^3(a, b, h)$ that have the same invariants of their Liouville foliation.

\section{Topology of isoenergy surfaces for Kovalevskaya case on $so(4)$}

In this section we determine the topological type of isoenergy surfaces $Q^3_{a, b, h}$ for the Kovalevskaya case on $so(4)$. Usually, for systems on $e(3)$, in order to solve this problem one can use the projection of an orbit $M^4 \approx T^*S^2$ on the Poisson sphere $S^2$ (see  [5] for more details). But for $so(4)$ the regular orbits are diffeomorphic to $S^2 \times S^2$ instead of $T^*S^2$ and we have to use another technique.

Here we show how the topology of $Q^3_{a, b, h}$ can be determined using the Fomenko's theory (see  [5]). The Fomenko-Ziechang invariants (also known as labeled molecules) for the Kovalevskaya case on $so(4)$ have been found in [8]. We will use notations from [8] in the statement and in the proof of the following theorem\footnote{Придётся это сказать, если не вставить молекулы.}.

\textbf{Theorem 1.} \textit{Any regular isoenergy manifold $Q^3_{a, b, h}$ for the integrable Kovalevskaya case on the Lie algebra $so(4)$ with Hamiltonian (1.3) and first integral (1.4) is homeomorhic to one of the following manifolds: \[S^3,\, 2S^3,\, S^1 \times S^2,\, 2 (S^1 \times S^2),\, \mathbb{R}P^3,\, (S^1 \times S^2) \# (S^1 \times S^2),\, (S^1 \times S^2) \# (S^1 \times S^2) \# (S^1 \times S^2).\] Moreover, for each labeled molecule 1-32 from [8] the topological type of the corresponding $Q^3_{a, b, h}$ is as in  Table~1.}

\begin{table}[h]
\centering
\begin{tabular}[t]{|c||l|c|l|}
\hline
Top. type of $Q^3$ & Fomenko-Ziechang invariants \\
\hline
  $S^3 $ & $1, 2, 3,  8, 9, 10, 11, 16, 17, 18, 21, 23$ \\
\hline
  $2 S^3 $ & $ 6, 7, 22 $ \\
\hline
  $ S^1 \times S^2 $ & $ 5, 12, 13, 14, 15, 20, 29, 32 $ \\
\hline
  $ 2 (S^1 \times S^2) $ & $ 27, 31 $ \\
\hline
  $\mathbb{R}P^3 $ & $ 4, 19, 24 $ \\
\hline
  $ (S^1 \times S^2) \# (S^1 \times S^2) $ & $ 26, 28, 30 $ \\
\hline
  $ (S^1 \times S^2) \# (S^1 \times S^2) \# (S^1 \times S^2) $ & $ 25 $ \\
\hline
\end{tabular}
\caption{Topology of  $Q^3_{a, b, h}$ for the Kovalevskaya case on $so(4)$.}
\end{table}

\textbf{Proof. } The proof is in several steps.

 1. First of all, let us note that for some labeled molecules the topology of $Q^3$ is well-known (see [5]). For instance, for the molecules $A \cfrac{r=0}{} A$ and $A \cfrac{r=\infty}{} A$ the surface $Q^3$ is diffeomorhic to $S^3$ and $S^1 \times S^2$ respectively. Thus we determine the topological type of molecules   $1, 7, 11, 13, 31$ from  [8].

Then, molecules $4, 28, 32$ were realized in the classical Kovalevskaya case and were denoted in [6] by $D, I, F$. Their isoenergy manifolds $Q^3 (D)$, $Q^3 (I)$  and $Q^3 (F)$ are diffeomorphic to $\mathbb{R}P^3, (S^1 \times S^2) \# (S^1 \times S^2)$ and $S^1 \times S^2$   respectively.

 2. Then, the topological type of an isoenergy surface $\left\{ h= \operatorname{const} \right\}$ does not change if we vary the value $H=h$ and do not pass through a critical point of $H$. Therefore we can identify  molecules by the following scheme.  Here we specify what molecules are identified for one of the bifurcation diagrams I-XI from [4] when passing through a point $y_2, y_5, y_6, y_8, y_9, y_{12}, z_2, z_8$. All the notations are as in [4] or [8].

\[1 \xrightarrow{y_2, V.8} 2 \xrightarrow{y_6, V.8} 3,\,\, 1 \xrightarrow{y_2, V.7} 16 \xrightarrow{y_9, III.4} 23; \qquad 22 \xrightarrow{y_8, III.1} 6 \xrightarrow{z_2, V.8} 7\] 

\[9 \xrightarrow{y_2, VII.4} 8 \xrightarrow{z_8, VII.4} 10,\,\, 9 \xrightarrow{z_8, VII.1} 11 \xrightarrow{y_2, VII.1} 10;\]

\[13 \xrightarrow{y_6, VII.6} 12 \xrightarrow{y_2, VII.6} 14,\,\, 13 \xrightarrow{y_2, VII.7} 15 \xrightarrow{y_9, III.1} 20;\]

\[17 \xrightarrow{y_5, IV.4} 18, \,\, 21 \xrightarrow{y_8, III.5} 18; \qquad 4 \xrightarrow{y_5, V.5} 19, \,\, 24 \xrightarrow{y_8, III.6} 19;\]

\[30 \xrightarrow{y_{13}, I.5} 26; \qquad 31 \xrightarrow{y_{13}, I.9} 27, \qquad 32 \xrightarrow{y_{13}, I.6} 29.\]

Thus we have determined the topology of all the molecules but  17, 18, 21 (these three molecules have the same topological type) and 25.

 3. Actually the Liouville equivalence can be defined not only for isoenergy manifolds of dynamic systems but on the set of foliated oriented manifolds with several directed circles where the fiber is not a smooth 2-manifold (e.g. the fiber is 1-dimensional or is like the special fiber of a 3-atom). Let us call these foliations as \textit{appropriate}.

Manifold with Liouville foliation of an integrable system can be described in terms of invariants of appropriate foliations that can be established on the manifold with such initial foliation. Cutting manifold and changing initial foliation on the produced pieces often  help to produce an appropriate foliation with simple molecule and well-known diffeomorphism type.

Neighbourhood of every singular fiber (i.e. a 3-atom) has structure of Seifert fibration with 2-dimensional base $M^2$ (foliated 2-manifold with boundaries, i.e. 2-atom). It is a direct product of $M^2$ and $S^1$ or such product factorized by involution $\mathbb{Z}_2$ on foliated $M^2$. Considered in our paper fibrations have type of direct product.

The following transformation simplifies foliation with every saddle atom of direct product type and will be explained in the case of atom $B$. Let us recall that a saddle atom $B$ has type of direct product and transforms two tori into one. Its base $M^2$ is a sphere with $3$ holes and foliation on $M^2$ has the only critical point of saddle type.

Elliptic atom $A$ is diffeomorphic to a solid torus foliated to a family of coaxial tori and one critical circle. It has a structure of direct product of $S^1$ and a foliated disc with one center critical point.

 Let $A$ and $B$ be glued s.th. the fiber of $B$ and the contractible cycle of $A$ are homologically equivalent or opposite (i.e. $r = \infty, \varepsilon = \pm 1$ on the edge $A - B$ of Fomenko--Zieschang invariant). The resulting foliation with base of $M^2$ will has fiber $S^1$ that is contracted to one point upon one of boundary circles.

 If neighbourhood of the saddle point contains two arcs of this circle then the base can be cut through the critical point s.th. the pullback of the cutting interval is a sphere $S^2$, see fig. 1.

 It allows to express initial manifold as a connected sum of two same manifolds. This idea was successfully used in [5] to the prove theorem 9.2 and proposition 4.5, e.g. to describe diffeomorphism class of $Q^3$ for a loop molecule of center-saddle type.

\begin{figure}[!htb]
\begin{center}
\minipage{0.70\textwidth}
\includegraphics[width=\linewidth]{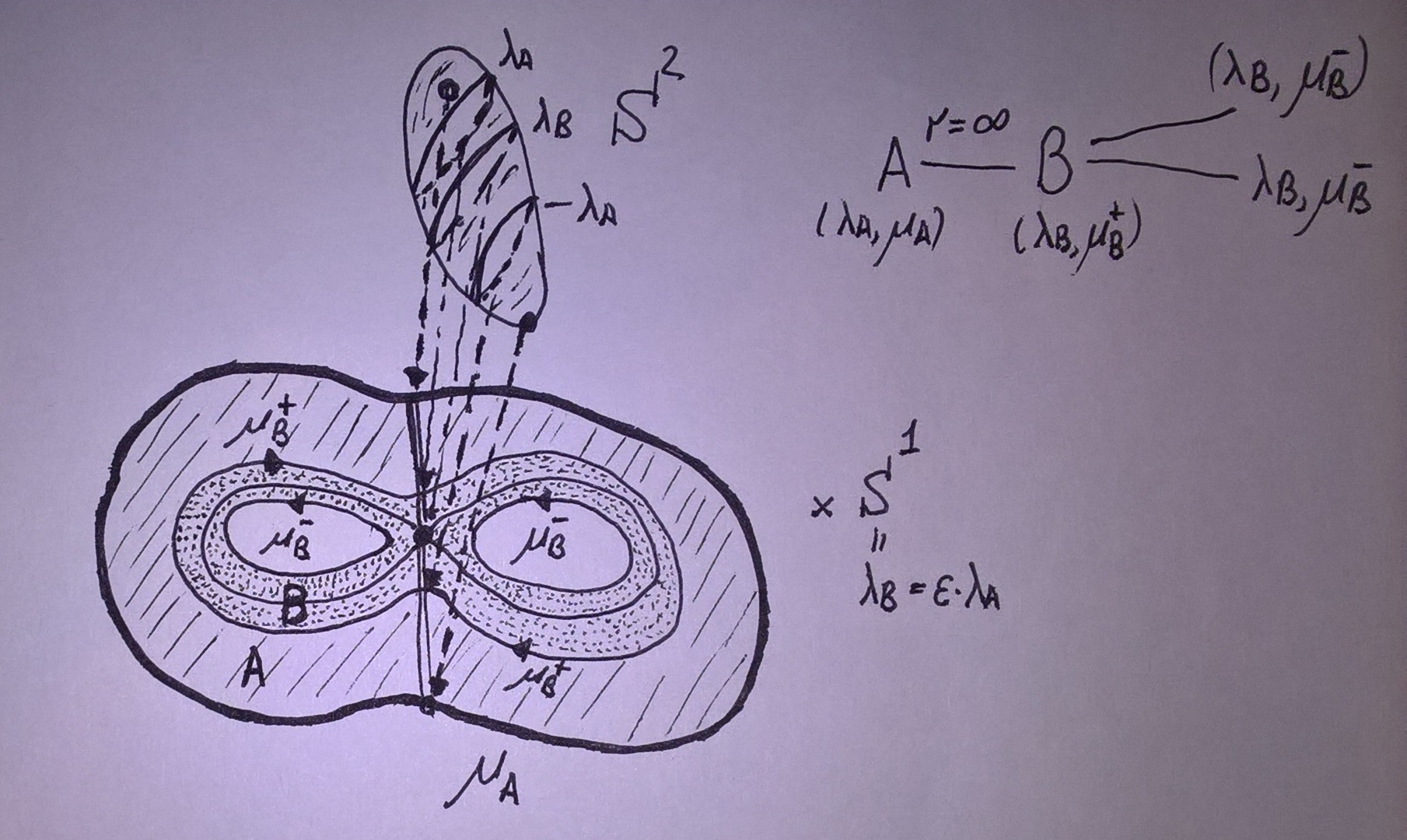}
   \caption{Preimage of an edge $A - B$ with $r = \infty$ as a connected sum}
\endminipage
\end{center}
\end{figure}

 To determine their diffeomorphism class one can preserve Liouville foliation on levels $K < k_{B}$ and establish foliation of $3$-atom $A$ on the rest part of each piece glued with 3-disk $D^3$. It means that set of degeneray of initial foliation should be the fiber of the new one and all regular fibers in the base of the new one should be closed, see fig. 2.

 The resulting foliation will have Fomenko--Zieschang invariant $A \,\, \cfrac{r = 0}{\varepsilon = 1} \,\, A^{*} \,\, \cfrac{r = 0}{\varepsilon =\pm 1}$ with $n = -1$ and thus is diffeomorphic to the isoenergy manifold of Goryachev- Chaplygin- Sretenskii integrable case, i.e. is class of $S^3$.

\begin{figure}[!htb]
\begin{center}
\minipage{0.5\textwidth}
\includegraphics[width=\linewidth]{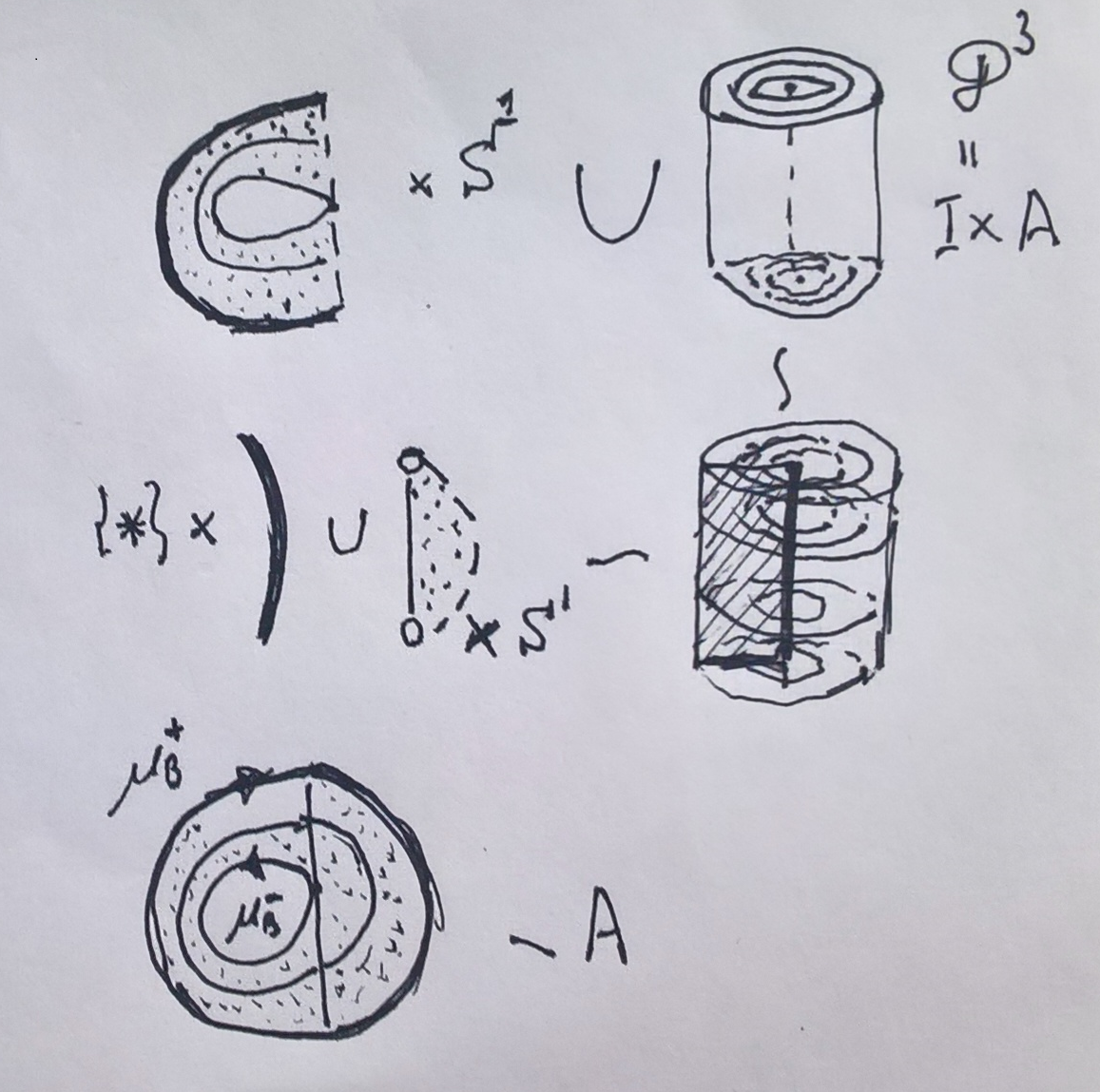}
   \caption{Introducing Liouville foliation on $Q^1 \# D^3$}
\endminipage
\end{center}
\end{figure}

 4. Let us consider the 3-atom $B$ connected with two other saddle atoms $B$ in the molecule $25$. Its mark $n$ is equal to $0$ and $r_i = 0$ for all three edges incidental to this $B$. It means that $\lambda$-cycles for three atoms incidental to this one correctly determine a section on this $B$ (there exists a section on 3-atom $B$ s.th. three $\lambda$-cycles are its intersections with boundary tori of $B$).

If a molecule has a fragment $A \cfrac{r = 0}{} V$ then the base of foliation can be considered as $S^2$ with two holes. The third hole is compensated by the base of $A$ atom. It will change the type of saddle atom. After some transformation of the foliation on the base (the gluing of two halves of the special fiber, fig.~3) one will have a molecule with two atoms $B$.

\begin{figure}[!htb]
\begin{center}
\minipage{0.80\textwidth}
\includegraphics[width=\linewidth]{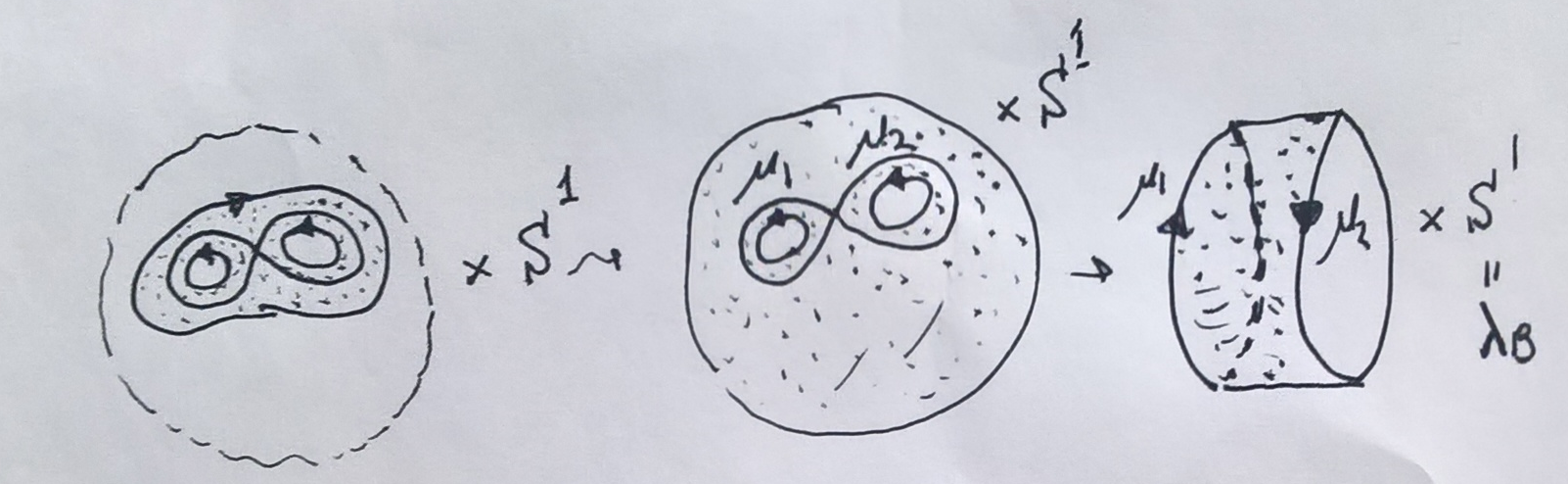}
   \caption{Transformation of foliation from $A - B$ with $r = 0$ to a critical torus}
\endminipage
\end{center}
\end{figure}

Topological type of $Q^3$ was not changed. The mark $r = \infty$ on the edge connecting them because $n$ was equal to $0$ on the disappeared $B$. It means that $\lambda$-cycles of the rest saddle atoms are homological to the fiber of disappeared $B$. The mark $\epsilon = -1$ because two $\mu$-cycles of atom $B$ should be oriented oppositely.

\begin{figure}[!htb]
\begin{center}
\minipage{0.65\textwidth}
\includegraphics[width=\linewidth]{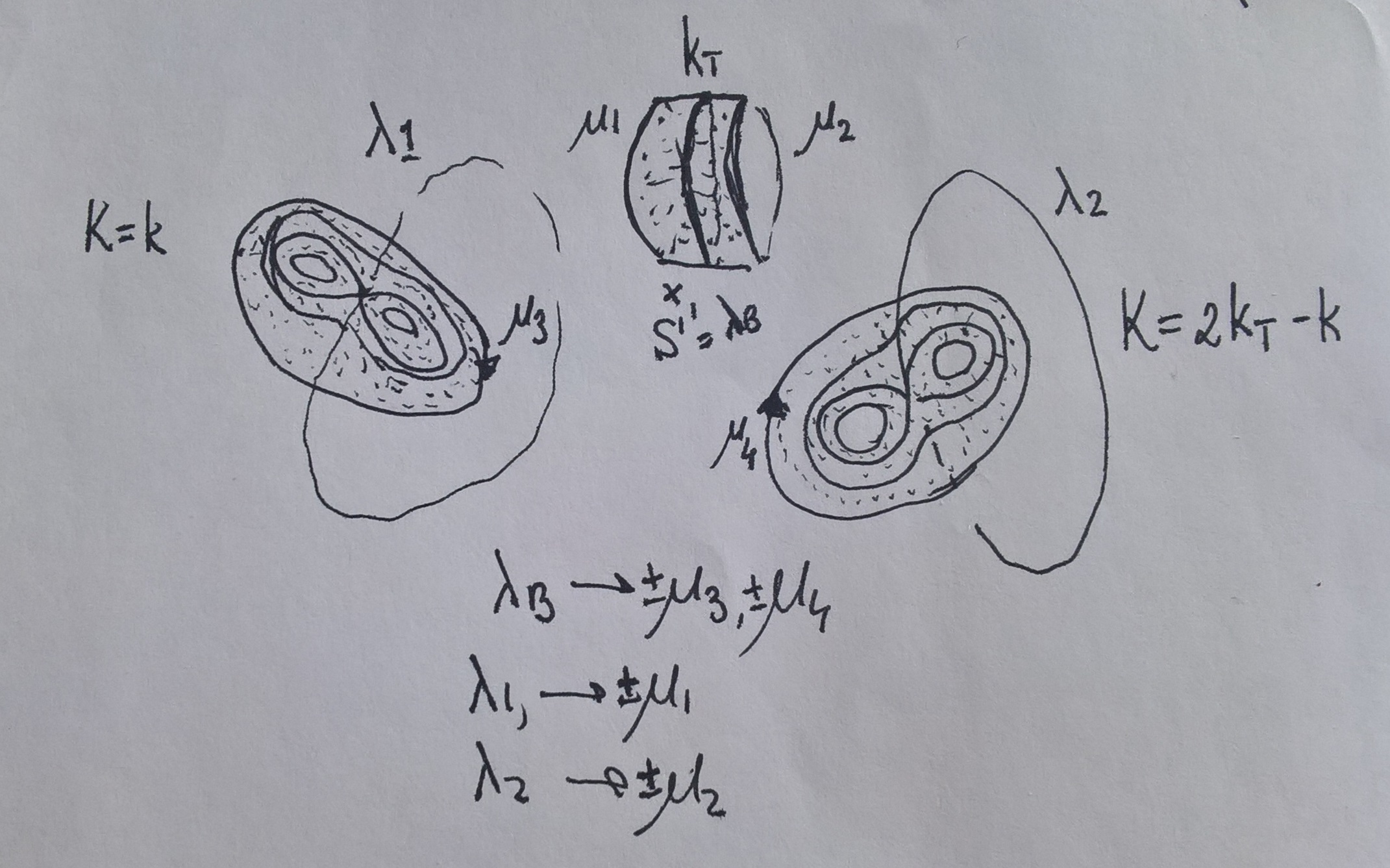}
   \caption{New foliation on $Q^3$ for the molecule 25}
\endminipage
\end{center}
\end{figure}

Two atoms $B \cfrac{r = \infty}{} B$ is a perturbation of atom $C_2$. The final molecule has for edges $A \cfrac{r = \infty}{} C_2$ and satisfies to theorem 9.2 in [5]. It was proved that the loop manifold (the preimage of a loop admissible curve) of center-saddle singularity $A \times C_2$ is diffeomorphic to a connected sum of $n+1$ numbers of $S^1 \times S^2$. Here $s = 2$ is the number of equilibria (points of the rank 0) on the special fiber of 2-atom $V = C_2$. It follows that $Q^3_{25}$ has the type of $(S^1 \times S^2) \# (S^1 \times S^2) \# (S^1 \times S^2)$.

\begin{figure}[!htb]
\begin{center}
\minipage{0.55\textwidth}
\includegraphics[width=\linewidth]{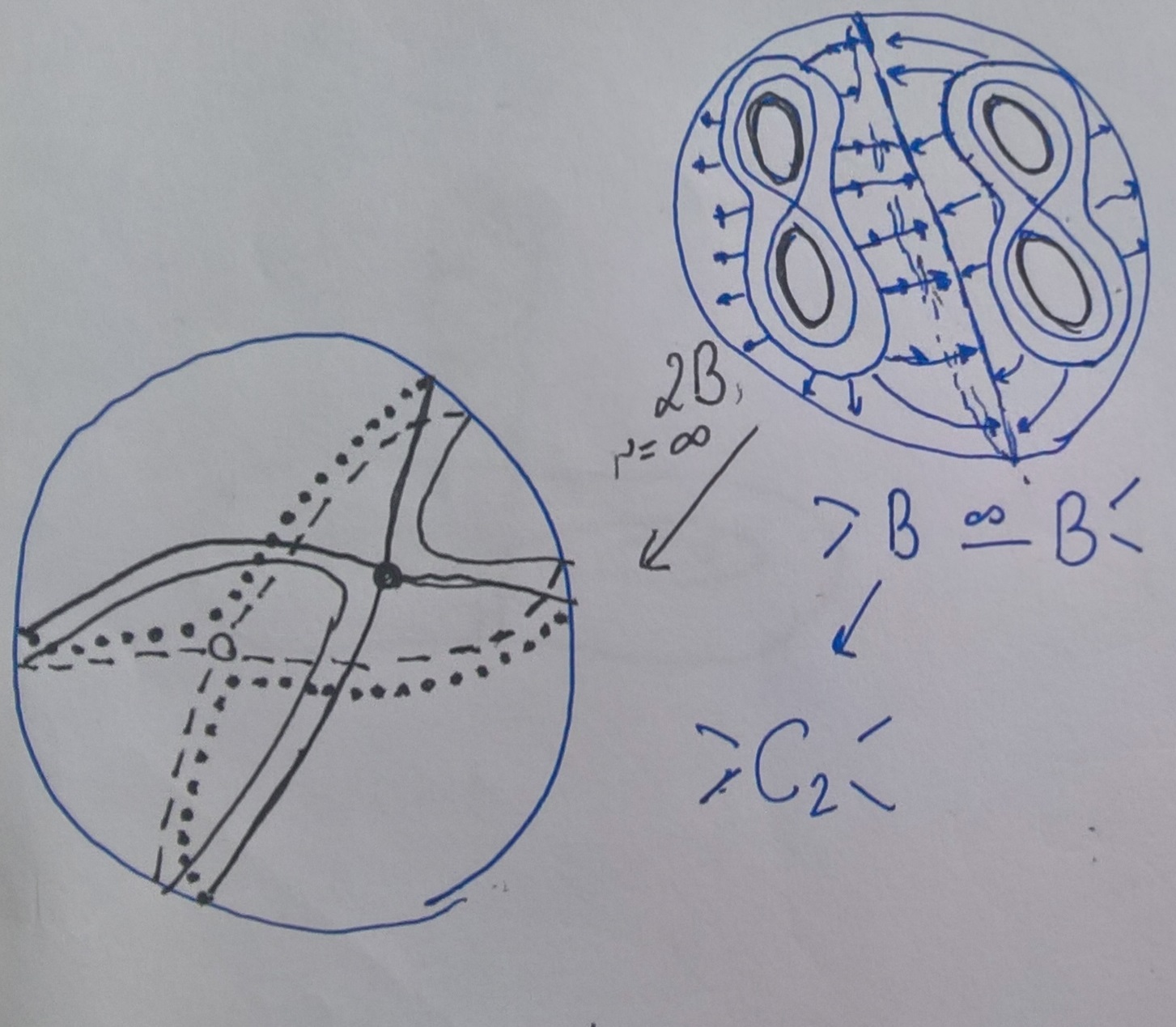}
   \caption{$B - B$ as a perturbation of $C_2$}
\endminipage
\end{center}
\end{figure}

\section{Acknowledgements}
This work was supported by Russian Scientific Foundation (project no. 17-11-01303). Author is grateful to I. Kozlov for fruitful discussions.

\end{fulltext}

\end{document}